\documentclass[10pt]{article}
\usepackage{amsmath, amssymb, graphicx, caption, float, algorithmic, algorithm, cases}
\usepackage{latexsym,amsfonts, subfigure, epstopdf, xcolor}

\numberwithin{equation}{section}
\numberwithin{figure}{section}
\numberwithin{table}{section}

\newtheorem{thm}{Theorem}[section]
\newtheorem{cor}{Corollary}[section]
\newtheorem{proposition}{Proposition}[section]
\newtheorem{lemma}{Lemma}[section]
\newtheorem{remark}{Remark}[section]

\newtheorem{assumption}{Assumption}[section]
\newcounter{nextauthor}
\begin{document}
\title{\Large {\bf Continuous and discrete-time accelerated methods for an
 inequality constrained convex optimization problem
\thanks{This work was supported by the National Natural Science Foundation of China (12171339, 12471296).
}}}
\author{Juan Liu$^a$, Nan-Jing Huang$^a$\thanks{Corresponding author: nanjinghuang@hotmail.com; njhuang@scu.edu.cn},
Xian-Jun Long$^b$ and Xue-song Li$^a$\\
{\small\it a. Department of Mathematics, Sichuan University, Chengdu, Sichuan 610064, P.R. China }\\
{\small\it b. College of Mathematics and Statistics, Chongqing Technology and Business University, Chongqing 400067, P.R. China}}
\date{}
\maketitle \vspace*{0mm}
\begin{center}
\begin{minipage}{5.5in}
{\bf Abstract.}
This paper is devoted to the study of acceleration methods for an inequality constrained convex optimization problem by using Lyapunov functions. We first approximate such a problem as an unconstrained optimization problem by employing the logarithmic barrier function. Using the Hamiltonian principle, we propose a continuous-time dynamical system associated with a Bregman Lagrangian for solving the unconstrained optimization problem. Under certain conditions,
we demonstrate that this continuous-time dynamical system exponentially converges to the optimal solution of the inequality constrained convex optimization problem. Moreover, we derive several discrete-time algorithms from this continuous-time framework and obtain their optimal convergence rates. Finally, we present numerical experiments to validate the effectiveness of the proposed algorithms.
\\ \ \\
{\bf Key Words:}
Inequality constrained convex optimization problem; Lyapunov function; Continuous-time dynamical system;
Discrete-time algorithm.
\\ \ \\
\end{minipage}
\end{center}

\section{Introduction}
Optimization problems naturally arise in various fields such as statistical machine learning, data analysis, economics, engineering, and computer science. Since Newton, numerous optimization methods have been proposed to solve minimization optimization problems,
including Newton methods, gradient descent, interior point methods, conjugate gradient methods and trust region methods (see, for example, Beck \cite{Ba}, Boyd et al. \cite{BPCPE,BV}, Nocedal and Wright \cite{NW}, Bertsekas et al. \cite{BNO},  Shor et al. \cite{S},  Esmaeili and Kimiaei \cite{HM}, Chen et al. \cite{CLO} and the references therein).
The rapid growth in the scale and complexity of modern datasets has led to a focus on gradient-based methods and also on the class of accelerated methods.
In 1983, Nesterov \cite{N} introduced the idea of acceleration in the context of gradient descent and showed that it achieve faster convergence rates than gradient descent.
Subsequently, the acceleration idea was applied to various optimization problems,
including nonconvex optimization \cite{ZL,KM}, non-Euclidean optimization \cite{N2,KBB}, quadratic programming \cite{LLM}, composite optimization \cite{BT,L,LLM1}, and stochastic optimization \cite{GL2,GL}.

To better understand the phenomenon of acceleration, researchers have provided insights into the design principles of accelerated methods by 
formulating the problem in continuous time and deriving algorithms through discretization.
Some researchers derive ordinary differential equations (ODEs) using  a limiting scheme based on existing algorithms.
Specifically,  Su et al. \cite{SBC} derived ODEs by taking continuous-time limits of existing algorithms and analyzed properties of the ODEs using the rich toolbox associated with ODEs, such as Lyapunov functions. Krichene et al. \cite{KBB} extended the results obtained in \cite{SBC} to non-Euclidean settings. Shi et al. \cite{SDJS} investigated an alternative limiting process that yields high-resolution differential equations.
Alternatively, 
some researchers adopt a variational point of view, deriving ODEs from a Lagrangian framework rather than from a limiting argument.
For instance,
inspired by the continuous-time analysis of mirror descent \cite{NY},
Wibisono et al. \cite{WWJ} derived ODEs directly from the underlying Lagrangian and established their stability using Lyapunov functions.
Further,
Jordan et al. \cite{BJ,MJ} proposed a systematic method for converting the ODEs presented in \cite{WWJ} into discrete-time algorithms by employing time-varying Hamiltonians and symplectic integrators.
Wilson et al. \cite{WRJ} established the equivalence between the estimate sequence technique and a family of Lyapunov functions in both continuous and discrete time, and used this connection to provide a simple and unified analysis of existing acceleration algorithms.

Notably, in either case, the continuous-time perspective provides the analytical power and intuition for the acceleration phenomenon as well as  design tools for developing new accelerated algorithms.
Some new acceleration methods for solving unconstrained optimization problems via the numerical discretization of ODEs have been obtained by, e.g., Chen et al. \cite{CSY}, Wang et al. \cite{WJW}, Luo and Chen \cite{LC,CL}, Bao et al. \cite{BCL,BCLS}, just mention a few.
A number of other papers have also contributed to the study of equality constraint optimization problems by working in continuous-time formulas.
For example,
Fazlyab et al. \cite{FKRP} extended the connection between acceleration methods and Lagrangian mechanics, as proposed in \cite{WWJ,WRJ}, to dual methods for linearly constrained convex optimization problems, and achieved exponential convergence in continuous-time dynamical system and an $O(1/k^2)$ convergence rate for discrete-time algorithms.
Zeng et al. \cite{ZLC} extended the continuous-time model from \cite{SBC} to equality constraint optimization problems, obtained a second-order differential system with primitive and multiplier variables, and proved its convergence rate $O(1/t^2)$. He et al. \cite{HHF2} and Attouch et al. \cite{ACFR} extended  the framework of \cite{ZLC} to separable convex optimization problems and discussed the convergence rate under different parameters.

As is well known, in optimization problems, inequality constraints allow the decision variables to vary within a certain range, rather than being strictly fixed to a particular value. This flexibility increases the difficulty of modeling and solving such problems. The discrete-time algorithms for solving inequality constrained optimization problems are discussed in some literature; see \cite{PTH,JPT,SR} and the references therein. However, to the best of our knowledge, there is no literature addressing  the convergence rate of accelerated methods for inequality constrained optimization problems from a continuous-time perspective. Therefore, it would be interesting and important to make an effort in this new direction.

The present paper is thus devoted to the study of new acceleration methods for solving an inequality constrained convex optimization problem from a continuous-time perspective. More precisely, we consider the following inequality constrained convex optimization problem (ICCOP):
\begin{align}
\min_{x\in\mathcal{X}} f(x), \;\; \mbox{s.t.}\;g_i(x)\leq 0, \;\forall i\in I. \label{1}
\end{align}
where $\mathcal{X}\subseteq \mathbb{R}^n$ is a convex compact set, $f$, $g_i:\mathcal{X}\rightarrow \mathbb{R}$ are real-valued convex functions with $i\in I:=\{1,\cdots,m\}$.
Compared with the known work addressing the equality constrained convex optimization problems via the dynamical systems, the primary challenge here is to conquer the trouble inequality constraints which requires certain novel approaches. Clearly, the methods used in those aforementioned papers do not apply here in a straightforward manner. Instead, one needs to carefully treat the inequality constraints, towards the acceleration methods. By notably utilizing barrier method and Lyapunov theory as well as the inequality techniques, we are able to propose some new acceleration methods for solving ICCOP. Our contributions in this paper are as follows:
\begin{itemize}
\item We approximate ICCOP as an unconstrained optimization problem by utilizing logarithmic barrier functions \cite[Sec.11]{BV}. This method differs from the dual reformulation of the linearly constrained convex problem in \cite{FKRP}, which merely yields a suboptimal gap.
\item We introduce a new Bregman-Lagrangian framework to derive a continuous-time dynamical system,
      where the potential energy is represented by an unconstrained optimization problem with a logarithmic barrier function, and the kinetic energy is defined by the Bregman divergence.
      Our Bregman-Lagrangian framework differs from the corresponding one in \cite{WWJ,WRJ,FKRP}.
\item Under certain conditions, we prove that the continuous-time dynamical system exponentially converges to the optimal solution.  Furthermore, we obtain an accelerated gradient method that matches the convergence rate of the underlying dynamical system by employing an implicit discretization method and propose another accelerated gradient method with a convergence rate of $O(1/k)$ by using hybrid Euler discretization that incorporates an additional sequence.
\end{itemize}

The rest of the paper is organized as follows.
In section 2,
we review relevant definitions and state the problem under consideration.
In section 3,
we define a Bregman Lagrangian on the continuous-time curves,
propose a continuous-time dynamical system,
and analysis its convergence properties.
In section 4,
we discretize the continuous-time dynamic to derive several discrete-time accelerated algorithms and investigate their convergence behavior.
In section 5,
we provide numerical experiments that validate the effectiveness of the proposed algorithms.
Finally,
conclusions are presented in section 6.

\section{Preliminaries}\label{s:main}

We first recall some notations and concepts. Let $\mathbb{R}$ and $\mathbb{R}_+$ be the set of real and nonnegative numbers,
respectively.
For nonnegative integer $k$, we use $x_k$ to denote discrete-time sequence. For $t\in\mathbb{T}\subseteq\mathbb{R}_+$, we use $X_t$ to denote continuous-time curves.
The notation $\dot{X}_ t$ means a derivative with respect to time,
i.e.,
$\dot{X}_ t = d X_t/ dt$.
Given a function $f:\mathbb{R}^n \rightarrow \mathbb{R}$,
its effective domain is defined as
$\mbox{dom}(f)=\{x\in \mathbb{R}^n: f(x)<\infty$\}.
A differentiable function $f$  is said to be
$L_f$-smooth with parameter $L_f>0$ if,
for any $x, y \in \mbox{dom}(f)$,
\begin{align}
f(y) \leq f(x) + \langle\nabla f(x),y-x\rangle+ \frac{L_f}{2}\|y- x\|^2.
\label{2.1}
\end{align}
This property is also equivalent to $\nabla f$ being Lipschitz continuous with parameter $L_f$.
A function $f$ is said to be $m_f$-strongly convex with parameter $m_f>0$ if,
for any $x, y \in \mbox{dom}(f)$,
\begin{align}
\langle\nabla f(x),y-x\rangle + \frac{m_f}{2}\|y-x\|^2+ f(x)\leq f(y).  \label{2.2}
\end{align}

Let $\mathcal{X}\subseteq \mathbb{R}^n$ be nonempty.
A function $h:\mathcal{X}\rightarrow\mathbb{R}$  is referred to as a distance-generating function with modulus $\mu>0$
if $h$ is continuously differentiable and $\mu$-strongly convex on $\mathbb{R}^n$.
The Bregman divergence associated with $h$ is defined as
\begin{align}
V_h(x,y)=h(x)-h(y)-\langle\triangledown h(y),x-y\rangle,\;
\forall x,y\in \mathcal{X}.
\label{5}
\end{align}
It follows from the definition that $V_h(x,y)\geq \frac{\mu}{2}\|x-y\|_2^2$ for all $x,y\in\mathcal{X}.$
Moreover, if $h(\cdot)=\frac{1}{2}\|\cdot\|_2^2$,
then $V_h(x,y)$  reduces to the Euclidean distance squared,
i.e.,
$V_h(x,y)=\frac{1}{2}\|x-y\|_2^2$.

\begin{remark}\label{rem:2.1}
\normalfont
If $\mathcal{X}$ is convex and compact,
then there exists a constant $M_{h,\mathcal{X}}>0$ such that $M_{h,\mathcal{X}} =\sup_ {x,y\in \mathcal{X}} V_h(x, y).$
This is a natural assumption imposed in constrained optimization problems \cite{L}.
\end{remark}

The Lagrange function associated to problem (\ref{1}) is formulated by
\begin{align}
\mathcal{L}(x,\lambda)=f(x)+\sum_{i\in I}\lambda_ig_i(x),
\label{2}
\end{align}
where $\lambda\in \mathbb{R}_+^{m}$ is a Lagrange multiplier.
The dual function associated to problem (\ref{1}) is defined as
$$G(\lambda)=\min_{x\in\mathcal{X}}\mathcal{L}(x,\lambda).$$
The set of all optimal dual solution is denoted
by
$\Lambda:=\arg \max_{\lambda\in \mathbb{R}^{m}} G(\lambda)$.
From now on, we assume that the minimizer $x^*$  of (\ref{1}) is unique, the objective function $f$ is $m_f$-strongly convex,
the constraint functions $g_i$ are convex for all $i\in I$,
and the problem (\ref{1}) is strictly feasible,
i.e., there exists $x\in\mathcal{X}$ such that
$g_i(x)<0$ for all $i\in I$.
Clearly,
the Lagrangian function $\mathcal{L}$ is twice differentiable,
strongly convex in $x$,
and concave in $\lambda$.
Furthermore,
the Slater's condition \cite{BNO} is satisfied.
Therefore,
there exists $\lambda^*\in\Lambda$ such that
the KKT conditions corresponding to the solution $x^*$ for (\ref{1}) are given by
\begin{align*}
&\nabla f(x^*)+\sum_{i\in I}\lambda_i^*g_i(x^*)=0, \\
&\lambda_i^*\geq 0, \quad -g_i(x^*)\ge 0, \quad \lambda_i^*g_i(x^*)=0, \;\;\forall i\in I.
\end{align*}

\begin{figure}[H]
  \centering
  \includegraphics[width=0.6\textwidth]{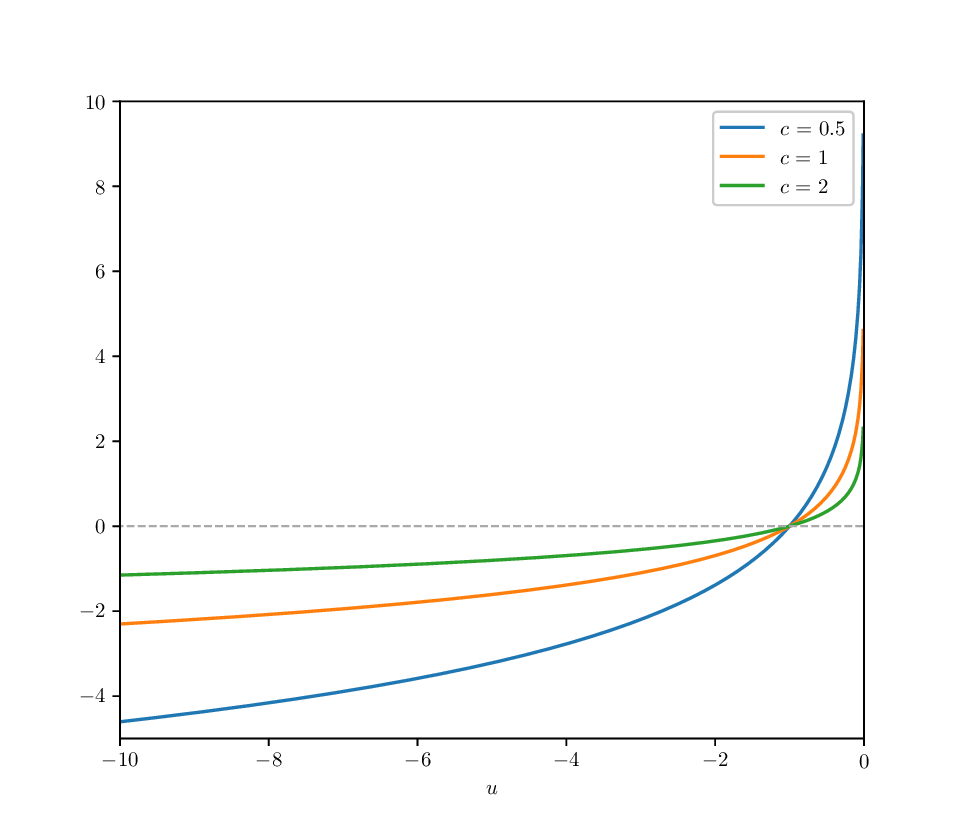}
  \caption{The dashed lines show the function $\mathbb{I}_-(u)$,
  and the solid curves show $\hat{\mathbb{I}}_-(u)= -(1/c) \log(-u)$, for $c = 0.5,\;1,\;2$.
  The curve for $c = 2$ gives the best approximation.}
  \label{fig:Figure 2.1}
\end{figure}

In what follows,
we approximate the constrained problem (\ref{1}) as an  unconstrained optimization problem using the logarithmic barrier function \cite[Sec. 11]{BV}.
Our first step is to rewrite the problem (\ref{1})
such that
the inequality constraints are implicit in the objective:
\begin{align}
\min_{x\in \mathcal{X}} f(x)+\sum_{i\in I}\mathbb{I}_-(g_i(x)),  \label{5}
\end{align}
where $\mathbb{I}_-:\mathbb{R}\rightarrow\mathbb{R}$ is the indicator function for the nonpositive reals, i.e.,
\begin{align*}
 \mathbb{I}_-(u)=\left\{
\begin{array}{ll}
0, &\quad\; u\leq0,\\
\infty,&\quad\; u>0.
\end{array}
\right.
\end{align*}
The problem (\ref{5}) has no inequality constraints, but its objective function is generally non-differentiable.
Consequently, we consider addressing this issue indirectly by introducing a differentiable approximation function.
Specifically, we approximate the indicator function $\mathbb{I}_-$ by the following function
$$\hat{\mathbb{I}}_-(u)=- \frac{1}{c}\sum_{i\in I}\log(-u),$$
where $c>0$ is a barrier parameter that controls the accuracy of the approximation.
Note that $\hat{\mathbb{I}}_-$ is convex, differentiable and increasing in $u$.
Figure \ref{fig:Figure 2.1} illustrates the function $\mathbb{I}_-$
and the approximation $\hat{\mathbb{I}}_-$,
for several values of $c$.
As $c$ increases,
the approximation becomes more accurate.
\setlength{\abovecaptionskip}{-2pt}

Substituting $\hat{\mathbb{I}}_-$ for $\mathbb{I}_-$ in (\ref{5}),
we obtain the following approximation:
\begin{align}
\min_{x\in K} f(x)- \frac{1}{c}\sum_{i\in I}\log(-g_i(x)),  \label{6}
\end{align}
where
$K:=\{x\in \mathcal{X}: g_i(x)<0, \forall i\in I \}$,
and the function $\phi(x)=(1/c)\sum_{i\in I}\log(-g_i(x))$ with $\mbox{dom}(\phi)=K$
is known as the logarithmic barrier for the problem (\ref{1}).
Since our goal is to design a dynamical system capable of indirectly tracking the optimal solution of the original problem (\ref{1}),
it is essential to initialize the dynamical system at a point inside $K$,
i.e., $X_0 \in K$.
To circumvent this restriction,
we introduce a nonnegative slack variable
$s$ into the approximation problem (\ref{6})
and reformulate it as follows:
\begin{align}
\min_{x\in \hat{K}} f(x)- \frac{1}{c }\sum_{i\in I}\log(s -g_i(x)),  \label{7}
\end{align}
where $\hat{K}:=\{x\in \mathcal{X}:g_i(x)<s ,\forall i\in I\}$ is an open and convex set that containing $K$.
For any $X_0 \in\mathcal{X}$,
we can choose $s> \max_{i\in I} g_i(X_0)$ to ensure $X_0\in \hat{K}$,
meaning the initial condition lies within the ``enlarged'' feasible set $\hat{K}$.
Let $\hat{x}^*$ be the minimizer of (\ref{7}).
To simplify the problem, we define the function
$$\Phi(x):= f(x)- \frac{1}{c }\sum_{i\in I}\log(s -g_i(x)).$$
Here, $\Phi$ is twice differentiable and strongly convex \cite{FPPR}.
Therefore,  the approximation problem (\ref{7}) is equivalent to
\begin{align}
\min_{x\in \hat{K}}\Phi(x).\label{8}
\end{align}
The optimal solution $\hat{x}^*$ of (\ref{7}) satisfies the optimality condition $\nabla\Phi(\hat{x}^* )=0$.

Next, we characterize the approximation error in terms of $c$, $s$,
and the optimal dual variable $\lambda^*$.

\begin{lemma}\label{lem:2.1}
\normalfont
Let $x^*$ be the minimizer of (\ref{1})
and $\hat{x}^*$ be the minimizer of  (\ref{8}).
Then,
for any $\lambda^*\in\Lambda$,
we have
\begin{align}
|f(\hat{x}^*)-f(x^*)|\leq \frac{m}{c }+s \|\lambda^*\|_1.
 \label{Lemma 1}
\end{align}
\end{lemma}
Proof:
Define $\tilde{x}^*$ as
\begin{align}
\tilde{x}^*:=\arg\min_{x\in \mathbb{R}^n} f(x)+\sum_{i\in I}\mathbb{I}_-(g_i(x)-s ),
\label{Lemma 1.1}
\end{align}
which is a perturbed version of the original optimization problem (\ref{1}) after including the slack variable $s $ in the constraints.
This problem coincides with the original problem (\ref{1}) when $s=0$.
If $s$ is positive,
it means that we have relaxed the $i$-th inequality constraint.
By perturbation and sensitivity analysis \cite[Sec.5.6]{BV},
we can establish the following inequality when $s \geq0$,
\begin{align}
0\leq f(x^*)-f(\tilde{x}^*)\leq \sum_{i\in I}\lambda_i^*s .
\label{Lemma 1.2}
\end{align}
The left inequality is based on the fact that the feasible set is enlarged when $s \geq0$,
and hence,
the optimal value does not increase.
The right inequality follows directly from a sensitivity analysis of the original problem \cite[Sec.5]{BV}.
On the other hand,
by replacing the indicator function $\mathbb{I}_-(u)$ in (\ref{Lemma 1.1}) with the log-barrier function $-(1/c) \log(-u)$,
we obtain the bound \cite[Sec.11]{BV}
\begin{align}
f(\hat{x}^*)-f(\tilde{x}^*)\leq \frac{m}{c }.
\label{Lemma 1.3}
\end{align}
Combining (\ref{Lemma 1.2}) and (\ref{Lemma 1.3}),
using triangle inequality,
we have
$$
|f(\hat{x}^*)-f(x^*)|\leq \frac{m}{c }+\sum_{i\in I}\lambda_i^*s =\frac{m}{c }+s \|\lambda^*\|_1.
$$
The proof is complete.\hfill
$\Box$

\begin{remark}\label{rem:2.2}
\normalfont
Lemma \ref{lem:2.1} shows that making suitable choices between $c$ and $s$ not only ensures that the right-hand side of inequation (\ref{Lemma 1}) converges to zero,
but also guarantees that the approximate solution $\hat{x}^*$ in (\ref{8}) converges to the optimal solution $x^*$ in (\ref{1}).
\end{remark}

\section{Continuous-time Dynamical System}\label{s:main}
In this section,
we design a continuous-time dynamical system and analyze its convergence properties.
We begin by deriving the dynamical system from the Bregman Lagrangian.
\subsection{The Bregman Lagrangian}\label{s:main}
Borrowing ideas from Wibisono et al. \cite{WWJ, WRJ},
define the Bregman Lagrangian as follows
\begin{align}
\mathbb{L}(X_t,\dot{X}_t,t)=e^{\alpha_t+\gamma_t}(V_h(X_t+e^{-\alpha_t}\dot{X_t},X_t)-e^{\beta_t}\Phi(X_t)),
\label{BL}
\end{align}
where $X_t\in \hat{K}$, $\dot{X}_t\in\mathbb{R}^n$,
and $t\in\mathbb{T}$ represent position, velocity and time,
respectively.
The functions
$\alpha,\gamma,\beta:\mathbb{T}\rightarrow\mathbb{R}$ are smooth increasing functions of time
that determine the weighting
of the velocity, the potential function, and the overall damping of the Lagrangian.
The ideal scaling conditions are given as follows

\begin{equation} \label{isc}
\left\{
\begin{aligned}
\dot{\beta}_t\leq e^{\alpha_t},\\
\dot{\gamma}_t=e^{\alpha_t},
\end{aligned}
\right.
\end{equation}
which are required in the stability analysis of continuous-time dynamical systems in the following paper.
For the Bregman Lagrangian $\mathbb{L}(X_t,\dot{X}_t,t)$,
denote by $J(\lambda)=\int_\mathbb{T}\mathbb{L}(X_t,\dot{X}_t,t)dt$  the action functional on curves $\{X_t:t\in\mathbb{T}\}$.
From Hamilton's principle \cite{BC},
finding the curve which minimizes the action functional is equivalent to finding a stationary point for the Euler-Lagrange equation,
that is
\begin{align}
\frac{d}{dt}  \frac{\partial \mathbb{L}}{\partial
\dot{X}_t}(X_t,\dot{X}_t,t) -\frac{\partial \mathbb{L}}{\partial
X_t}(X_t,\dot{X}_t,t)=0.
\label{EL}
\end{align}

\begin{proposition} \label{pro:3.1}
\normalfont
Under the ideal scaling conditions (\ref{isc}),
the Euler-Lagrange equation (\ref{EL}) associated with the Bregman Lagrangian (\ref{BL}) reduces to the following continuous-time dynamical system
\begin{align}
d\nabla h(X_t+e^{-\alpha_t}\dot{X}_t)=-e^{\alpha_t+\beta_t}\nabla \Phi(X_t)dt.
\label{sODE}
\end{align}
\end{proposition}
Proof.
The partial derivatives of the Bregman Lagrangian can be weitten,
\begin{align*}
\frac{\partial \mathbb{L}}{\partial \dot{X}_t}(X_t,\dot{X}_t,t)
&=
e^{\gamma_t}(Y_t-\nabla h(X_t)), \\
\frac{\partial \mathbb{L}}{\partial
X_t}(X_t,\dot{X}_t,t)
&=
e^{\alpha_t}\frac{\partial \mathbb{L}}{\partial \dot{X}_t}(X_t,\dot{X}_t,t)-
e^{\gamma_t}\dot{X}_t\frac{d}{dt} \nabla h(X_t)
-e^{\alpha_t+\beta_t+\gamma_t}
\nabla \Phi(X_t).
\end{align*}
We also compute the time derivative of $\frac{\partial \mathbb{L}}{\partial
\dot{X}_t}(X_t,\dot{X}_t,t)$,
\begin{align*}
\frac{d}{dt}
\frac{\partial \mathbb{L}}{\partial
\dot{X}_t}(X_t,\dot{X}_t,t)
&=
\dot{\gamma}_t
\frac{\partial \mathbb{L}}{\partial \dot{X}_t}(X_t,\dot{X}_t,t)+
e^{\gamma_t}
(\frac{d}{dt} Y_t-\dot{X}_t\frac{d}{dt} \nabla h(X_t)).
\end{align*}
Using the ideal scaling condition (3.2),
the Euler-Lagrange equation given by
$$d\nabla h(X_t+e^{-\alpha_t}\dot{X}_t)=-e^{\alpha_t+\beta_t}\nabla \Phi(X_t)dt.$$
In the calculation,
the terms involving
$\frac{d}{dt} \nabla h(X_t)$ are eliminated
and the terms involving
$\frac{\partial \mathbb{L}}{\partial
\dot{X}_t}(X_t,\dot{X}_t,t)$
are simplified.\hfill
$\Box$

\begin{remark}\label{rem:3.1}
\normalfont
For the dynamical system (\ref{sODE}),
the solution $X_t$ must ensure that the argument
$s-g_i(X_t)$ of the logarithmic barrier function (\ref{8}) remains positive.
That is,
for all $t\geq 0$, $X_t \in \hat{K}$.
In deed,
$$\nabla\Phi(X_t)=\nabla f(X_t)+\frac{1}{c }\sum_{i\in I}\frac{\nabla g_i(X_t)}{s -g_i(X_t)},$$
if exists $t\geq 0$
such that $s=g_i(X_t)$,
then $\|\nabla\Phi(X_t)\|_2$ is unbounded (singular) at the boundary of $\hat{K}$,
which is impossible since the strong convexity of $\Phi$ implies that
$\|\nabla\Phi(X_t)\|_2$ is bounded for all $t\geq 0$.
Thus, we must have that
$X_t \in \hat{K}$ for all $t\geq 0$.
\end{remark}

\subsection{Convergence rates of the Continuous-time Dynamical System}\label{s:main}

In this subsection,
we study the convergence of the continuous-time dynamical system given in Proposition \ref{pro:3.1}.
First,
we show that the solution of the continuous-time dynamical system (\ref{sODE}) exponentially converts to the approximate solution $\hat{x}^*$ in (\ref{7}) by employing Lyapunov function approach \cite{LF}.

\begin{lemma}\label{lem:3.2}
\normalfont
Let $X_t$ be the solution of (\ref{sODE}) and $\hat{x}^*$ be the minimizer of (\ref{8}).
Then, under the ideal scaling conditions (\ref{isc}), one has $\Phi(X_t)-\Phi(\hat{x}^*)= O( e^{-\beta_t})$.
In addition, for large enough $c$,
\begin{align}
 \|X_t- \hat{x}^*\|_2 \leq \sqrt{\frac{2}{m_f}\mathcal{E}_0} e^{-\frac{\beta_t}{2}},
 \label{TH1}
\end{align}
where
$\mathcal{E}_0=e^{\beta_0}(\Phi(X_0)-\Phi(\hat{x}^*))+V_{h}(\hat{x}^*,X_0+e^{-\alpha_0}\dot{X}_0)$.
\end{lemma}
Proof.
Consider the following Lyapunov function
\begin{align}
\mathcal{E}_t=e^{\beta_t}(\Phi(X_t)-\Phi(\hat{x}^*))+V_{h}(\hat{x}^*,X_t+e^{-\alpha_t}\dot{X}_t).
\label{LF}
\end{align}
Observe that (\ref{LF}) is non-negative.
since the function $h$ is convex,
we have $V_{h}(x, y) \geq 0$ for all $x, y \in \mathbb{R}^n$.
The rest of the energy functional is positive due to $\Phi(x)-\Phi(\hat{x}^*)\geq0$ for any $x\in \mathbb{R}^n$.
The time derivative of the Lyapunov functional is
\begin{align*}
\dot{\mathcal{E}}_t
=\dot{\beta}_t e^{\beta_t}(\Phi(X_t)-\Phi(\hat{x}^*))+
e^{\beta_t}\langle\Phi(X_t),\dot{X}_t\rangle+
e^{\alpha_t+\beta_t}\langle\nabla\Phi(X_t),\hat{x}^*-X_t-e^{-\alpha_t}\dot{X}_t\rangle.
\end{align*}
Substitute (\ref{sODE}) into the above formula,
we have
\begin{align*}
\dot{\mathcal{E}}_t
=-e^{\alpha_t+\beta_t}V_{\Phi}(\hat{x}^*,X_t)+
(\dot{\beta}_t-e^{\alpha_t})e^{\beta_t}(\Phi(X_t)-\Phi(\hat{x}^*)),
\end{align*}
where $V_{\Phi}(\hat{x}^*,X_t)=\Phi(\hat{x}^*)-\Phi(X_t)-\langle\nabla\Phi(X_t),\hat{x}^*-X_t\rangle$
is the Bregman divergence of $\Phi$.
By the ideal scaling conditions (\ref{isc}) and the nonnegativity of $V_{\Phi}$,
we have $\dot{\mathcal{E}}_t\leq0$.
It follows that
$\mathcal{E}_t\leq\mathcal{E}_0$
for all $t\geq t_0\in \mathbb{T}$.
This together with the definition of $\mathcal{E}_t$
and the fact of $V_{h}(\hat{x}^*,X_t+e^{-\alpha_t}\dot{X}_t) \geq 0$ yields
\begin{align}
\Phi(X_t)-\Phi(\hat{x}^*)\leq \mathcal{E}_0 e^{-\beta_t}=O( e^{-\beta_t}).
 \label{TH1-1}
\end{align}
Furthermore,
since $f$ is $m_f$-strongly convex,
and $c>0$ large enough,
it follows that $\Phi$ is $m_f$-strongly convex for $X_t\in\hat{K}$ and so
\begin{align}
 \frac{m_f}{2}\|X_t- \hat{x}^*\|_2^2 \leq \Phi(X_t) -\Phi(\hat{x}^*).
 \label{TH1-2}
\end{align}
Combining (\ref{TH1-1}) and (\ref{TH1-2}),  (\ref{TH1}) holds.
\hfill
$\Box$

In order to establish the convergence relationship between the solution generated by continuous-time dynamical system (\ref{sODE}) and the optimal solution in (\ref{1}),
it is necessary to make
the following assumptions about the barrier parameter $c$,
the slack variable $s$ and the optimal dual variable $\lambda^*$.

\begin{assumption}\label{a:3.1}
\normalfont
For $t\in \mathbb{T}$, assume that
\begin{itemize}
\item[(i)] $c$ is a time-dependent positive barrier function, i.e., $c(t):=e^{\beta_t}$;
\item[(ii)] $s$ is a nonnegative time-dependent slack function,  i.e., $s(t):=e^{-\beta_t}$;
\item[(iii)] the optimal dual variables satisfy $\lim_{t\rightarrow\infty}(\inf_{\lambda^*\in\Lambda}\|\lambda^*\|_1)e^{-\beta_t}=0$.
\end{itemize}
\end{assumption}

\begin{remark}\label{rem:3.2}
\normalfont
If Assumption \ref{a:3.1} holds, then one has the following results:  (i) the barrier function $c$ asymptotically diverges to infinity; (ii) the slack  function $s$ goes to zero exponentially fast; (iii) the approximation error in (\ref{Lemma 1}) vanishes asymptotically providing $\inf_{\lambda^*\in\Lambda}\|\lambda^*\|_1$ is bounded.
\end{remark}

Under Assumption \ref{a:3.1}, we obtain the following important theorem by utilizing Lemmas \ref{lem:2.1} and \ref{lem:3.2}.

\begin{thm}\label{thm:3.1}
\normalfont
Let $x^*$, $\hat{x}^*$ and $X_t$ be the solution of (\ref{1}), (\ref{8}) and (\ref{sODE}),
respectively. Then, under Assumption \ref{a:3.1}, one has $\lim_{t\rightarrow\infty}\|X_t-x^*\|_2=0.$
In addition, if $f$ is $L_f$-smooth, then
$$\|X_t-x^*\|_2 \leq C_0e^{-\frac{\beta_t}{2}}=O(e^{-\frac{\beta_t}{2}}),$$
where $C_0:=\sqrt{\frac{2}{m_f}\mathcal{E}_0} +\frac{2}{m_f}\sqrt{2L_f(m+ \|\lambda^*\|_1)}$.
\end{thm}
Proof: Under Assumption \ref{a:3.1}, by Lemmas \ref{lem:2.1} and \ref{lem:3.2},
we can directly get $\lim_{t\rightarrow\infty}\|X_t-x^*\|_2=0$.
Next,
we discuss the convergence rate.
By the strong convexity of $f$,
\begin{align*}
f(x^*)&
\geq f(\hat{x}^*)+ \langle\nabla f(\hat{x}^*),x^*-\hat{x}^*\rangle+ \frac{m_f}{2} \|x^*-\hat{x}^*\|_2^2\\
&\geq f(\hat{x}^*)-\|\nabla f(\hat{x}^*)\|_2 \|(x^*-\hat{x}^*)\|_2+ \frac{m_f}{2} \|x^*-\hat{x}^*\|_2^2.
\end{align*}
It follows from $f(x^*)\leq f(\hat{x}^*)$ that
\begin{align}
 \frac{m_f}{2} \|x^*-\hat{x}^*\|_2 \leq \|\nabla f(\hat{x}^*)\|_2.
 \label{TH3.1}
\end{align}
Since $f$ is $L_f$-smooth,
for any $x, y\in \mbox{dom}(f)$,
$f(x^*)
\leq f(y)
\leq f(x)+ \langle\nabla f(x),y-x\rangle+ \frac{L_f}{2} \|y-x\|_2^2$.
Let $x=\hat{x}^*$
and $y=-\frac{\nabla f(\hat{x}^*)}{L_f}+\hat{x}^*$,
we get
\begin{align}
\|\nabla f(\hat{x}^*)\|_2\leq \sqrt{2L_f(f(\hat{x}^*)-f(x^*) )}.
 \label{TH3.2}
\end{align}
Combining (\ref{TH3.1}) and (\ref{TH3.2}) yields
$$\|x^*-\hat{x}^*\|_2 \leq \frac{2}{m_f}\sqrt{2L_f(f(\hat{x}^*-f(x^*) )}.$$
By Assumption \ref{a:3.1} and Lemma \ref{lem:2.1}, one has
$|f(\hat{x}^*)-f(x^*)|\leq (m+ \|\lambda^*\|_1)e^{-\beta_t}$.
It follows that
$$0\leq\|\hat{x}^*-x^*\|_2\leq \frac{2}{m_f}e^{-\frac{\beta_t}{2}}\sqrt{2L_f(m+  \|\lambda^*\|_1)}.$$
By Lemma \ref{lem:3.2}, we have
$$0\leq\|X_t-x^*\|_2\leq\|X_t-\hat{x}^*\|_2+\|\hat{x}^*-x^*\|_2\leq C_0e^{-\frac{\beta_t}{2}}.$$
This completes the proof.\hfill
$\Box$
\begin{remark}\label{rem:3.3}
\normalfont
Theorem \ref{thm:3.1} shows that the solution of the continuous-time dynamical system (\ref{sODE}) asymptotically converges to the optimal solution in (\ref{1}).
For any positive constant $p$,
if $\beta_t = 2p \log t$,
then we attain a polynomial convergence rate $O(1/t^p)$,
if $\beta_t=2tp$,
then we attain an exponential convergence rate $O(e^{-tp})$.
\end{remark}

\section{Discrete-time Algorithms }\label{s:main}
Wibisono et al. \cite{WWJ} pointed out
that the numerical algorithms whose convergence rates match those of the underlying dynamical system cannot be easily obtained by using the naive numerical discretization method,
and even the proposed discretization algorithm fails to converge.
Therefore,
it is meaningful to study the discrete-time algorithms.

In this section,
we design several acceleration algorithms for the optimization problem (\ref{8}) based on the Euler discretization \cite{KP} of dynamical system (\ref{sODE})
and analyze the convergence of the proposed discrete-time algorithms
using discrete-time Lyapunov functions.
This ideal is inspired by Wilson et al. \cite{WRJ},
which discussed various discretization schemes for the continuous-time dynamical system of unconstrained optimization.
From now on,
we assume that the first ideal scaling condition in (\ref{isc}) holds with equality,
i.e., $\dot{\beta}_t=e^{\alpha_t}$.
In order to get the discretizations of continuous-time dynamical system (\ref{sODE}),
we first rewrite it as the following system of first-order equations
\begin{equation} \label{first order equtions}
\left\{
\begin{aligned}
&Z_t=X_t+\frac{1}{\dot{\beta}_t}\dot{X}_t,\\
&d\nabla h(Z_t)=-\dot{\beta}_te^{\beta_t}\nabla \Phi(X_t)dt.
\end{aligned}
\right.
\end{equation}
Next,
we discretize $X_t$ and $Z_t$ into sequences $x_k$ and $z_k$ with time step $\delta>0$.
That is,
we make the identification $t = \delta k$
and set
$$x_k=X_t,\;\;
x_{k+1}=X_{t+\delta},\;\;
z_k=Z_t,\;\;
z_{k+1}=Z_{t+\delta}.$$
The explicit (forward) Euler discretization for
$\dot{X}_t$ and $\frac{d}{dt}\nabla h(Z_t)$
are defined as
$$\frac{x_{k+1}-x_k}{\delta}\approx \dot{X}_t,\;
\frac {\nabla h(z_{k+1})-\nabla h(z_{k})}{\delta}\approx \frac{d}{dt}\nabla h(Z_t).$$
The implicit (backward) Euler discretization for
$\dot{X}_t$ and $\frac{d}{dt}\nabla h(Z_t)$
are defined as
$$\frac{x_{k}-x_{k-1}}{\delta}\approx\dot{X}_t,\;
\frac {\nabla h(z_{k})-\nabla h(z_{k-1})}{\delta}\approx\frac{d}{dt}\nabla h(Z_t).$$
For the scaling parameters,
we choose $A_k = e^{\beta_t}$,
and the discretizations as follows
$$\alpha_k:=\frac{A_{k+1}-A_{k}}{\delta}\approx\frac{d e^{\beta_t}}{dt},\;
\tau_k:=\frac{A_{k+1}-A_{k}}{\delta A_{k}}\approx\dot{\beta}_t=\frac{\frac{d}{dt}e^{\beta_t}}{e^{\beta_t}}.$$

\subsection{Implicit Euler Discretization}\label{s:main}
In this subsection,
we show that the implicit discretization of continuous-time dynamical system (\ref{first order equtions}) can produce an acceleration algorithm whose convergence rate matches that of the underlying dynamical system.

\begin{algorithm}
  \caption{Accelerated  Gradient Method (AGM1)}
  \label{algo:threeline}
  \begin{algorithmic}
  \STATE \;\textbf{Initialization:} Choose $x_0=z_0\in \mathbb{R}^n $, $\delta>0$.
  \STATE \;\textbf{For} $k=0,1,\cdots$ \textbf{do}
  \STATE \; \quad  \textbf{Step 1:} Choose $A_k$. Set
          \STATE \; \quad \quad\quad \quad\quad
          $c_k=A_k$, $s_k=1/A_k$, $\alpha_k=\frac{A_{k+1}-A_{k}}{\delta}$,
           $\tau_k=\frac{A_{k+1}-A_{k}}{\delta A_{k}}$.
    \STATE \; \quad \textbf{Step 2:} Compute
    $x_{k+1}=\frac{\delta\tau_k}{1+\delta\tau_k} z_{k+1}+ \frac{1}{1+\delta\tau_k}x_{k}$.
    \STATE \;   \quad \textbf{Step 3:} Compute
     $z_{k+1}=\arg\min_{z\in \hat{K}}\{\alpha_k\langle \nabla\Phi(x_{k+1}),z\rangle+ \frac{1}{\delta}V_h(z,z_k)\}$.
     \STATE \; \quad  \textbf{If} A stopping condition is satisfied  \textbf{then}
      \STATE \; \quad \quad \quad\textbf{Return} $(x_k,z_k) $
       \STATE \;\quad   \textbf{end}
    \STATE  \textbf{end}
  \end{algorithmic}
\end{algorithm}

By the implicit Euler discretization of $\dot{X}_t$ and $\frac{d}{dt}\nabla h(Z_t)$,
the continuous-time dynamical system (\ref{first order equtions})  can be discretized as follows
\begin{equation} \label{discretized.1}
\left\{
\begin{aligned}
&z_{k+1}=x_{k+1}+ \frac{1}{\delta\tau_k}(x_{k+1}-x_{k}), \\
&\nabla h(z_{k+1})-\nabla h(z_{k})=-(A_{k+1}-A_k)\nabla\Phi (x_{k+1}).
\end{aligned}
\right.
\end{equation}
It is worth noting that (\ref{discretized.1}) actually gives the optimality condition for a discrete-time algorithm,
which is depicted in Algorithm 1.

The following theorem shows that the convergence behavior of Algorithm 1 by using discrete-time Lyapunov function.

\begin{thm}\label{thm:4.1}
\normalfont
Suppose $f$ is $L_f$-smooth.
Define the discrete-time Lyapunov function
\begin{align}
\mathcal{E}_k=A_k(\Phi(x_k)-\Phi(\hat{x}^*))+V_{h}(\hat{x}^*,z_k).
\label{LF2}
\end{align}
Then,
under Assumption \ref{a:3.1},
$\mathcal{E}_{k+1}-\mathcal{E}_{k}\leq0$
and
$$\|x_k-x^*\|_2=O(\frac{1}{\sqrt{A_k}}).$$
\end{thm}
Proof.
From (\ref{LF2}),
we have
\begin{align*}
 \mathcal{E}_{k+1}-\mathcal{E}_{k}
& =  A_{k+1}(\Phi(x_{k+1})-\Phi(\hat{x}^*))+V_{h}(\hat{x}^*,z_{k+1})-A_k(\Phi(x_k)-\Phi(\hat{x}^*))-V_{h}(\hat{x}^*,z_k)\\
& = h(\hat{x}^*)-h(z_{k+1})-\langle\nabla h(z_{k+1}),\hat{x}^*-z_{k+1}\rangle-h(\hat{x}^*)+h(z_{k})+\langle\nabla h(z_{k}),\hat{x}^*-z_{k}\rangle\\
  & \;\;\;\;
  +(A_{k+1}-A_k)(\Phi(x_{k+1})-\Phi(\hat{x}^*))+A_k(\Phi(x_{k+1})-\Phi(x_k))\\
& = (A_{k+1}-A_k)\langle\nabla\Phi(x_{k+1}),\hat{x}^*-x_{k+1}\rangle+A_k\langle\nabla\Phi(x_{k+1}),x_k-x_{k+1}\rangle-V_{h}(z_{k+1},z_k) \\
  & \;\;\;\;
  +(A_{k+1}-A_k)(\Phi(x_{k+1})-\Phi(\hat{x}^*))+A_k(\Phi(x_{k+1})-\Phi(x_k))\\
& \leq 0.
\end{align*}
The second equality follows from the definition of $V_h$, the third equality follows from (\ref{discretized.1}),
and the inequality follows from the convexity of $\Phi$
and the fact of $V_h\geq 0$.
Summing from $k=0$ to $k-1$,
we obtain $\mathcal{E}_k\leq \mathcal{E}_0$.
By definition of $\mathcal{E}_k$, one has $\Phi(x_k)-\Phi(\hat{x}^*)\leq
\frac{\mathcal{E}_0}{A_k}$. This together with the strong convexity of $\Phi$ yields
$\|x_k-\hat{x}^*\|_2\leq \sqrt{\frac{2\mathcal{E}_0}{m_fA_k}}$.
For $A_k = e^{\beta_t}$, by Assumption \ref{a:3.1} and Lemma \ref{lem:2.1}, one has
$|f(\hat{x}^*)-f(x^*)|\leq \frac{m+ \|\lambda^*\|_1}{A_k}$.
Since $f$ is strongly convex and smooth,
using the same arguments as the proof of Theorem 3.1,
we obtain
$$\|x_k-x^*\|_2\leq \sqrt{\frac{2\mathcal{E}_0}{m_fA_k}}+ \frac{2}{m_f}\sqrt{\frac{2L_f(m+ \|\lambda^*\|_1)}{ A_k}}=O(\frac{1}{\sqrt{A_k}}).$$
This completes the proof.\hfill
$\Box$

\begin{remark}\label{rem:4.1}
\normalfont
Theorem \ref{thm:4.1} establishes an acceleration gradient algorithm that matches the convergence rate of the dynamic system (\ref{first order equtions}) by choosing an appropriate implicit discretization method.
For constant $p > 0$.
If $A_k = (\delta k)^{2p}$,
then achieve a polynomial convergence rate $O(\frac{1}{(\delta k)^p})$.
If $A_k =e^{2\delta kp}$ for some constant $p>0$,
then achieve an exponential convergence rate $O(e^{-\delta kp})$.
\end{remark}

\begin{remark}\label{rem:4.2}
\normalfont
If $g_i(x)=0$ for all $i\in I$ ,
then Algorithm 1 can solve unconstrained convex optimization problems and achieve the best convergence rate.
\end{remark}

\subsection{Hybrid Euler Discretization with an Additional Sequence}\label{s:main}
Although Algorithm 1 is closely related to the continuous-time dynamical system (\ref{first order equtions}) and has an optimal convergence rate,
it involves an implicit update step (Step 3) which is pretty hard to solve in practice.
Thus,
it's natural to consider whether computationally efficient algorithms can be derived using an explicit Euler discretization of one of the sequence.
In this subsection,
we present two algorithms using the aforementioned technique.
One directly yields a gradient method,
and the other,
with an additional sequence,
yields an accelerated gradient method that differs from the previous section.

First,
we show Algorithm 2,
which combines implicit and explicit Euler discretizations of dynamical system (\ref{first order equtions}) with an additionally updating sequence.
\begin{algorithm}
  \caption{Accelerated  Gradient Method (AGM2)}
  \label{algo:threeline}
  \begin{algorithmic}
  \STATE \;\textbf{Initialization:} Choose $y_0=z_0\in \mathbb{R}^n $, $\delta>0$, $C>0$, $\eta>0$, $A_0=C$.
  \STATE \;\textbf{For} $k=0,1,\cdots$ \textbf{do}
  \STATE \; \quad  \textbf{Step 1:} Set
          $A_{k+1}=C(k+1)^2$, $c_k=A_k$, $s_k=1/A_k$, $\alpha_k=\frac{A_{k+1}-A_{k}}{\delta}$,
           $\tau_k=\frac{A_{k+1}-A_{k}}{\delta A_{k}}$.
    \STATE \; \quad \textbf{Step 2:} Compute
    $x_{k+1}=\delta\tau_k z_{k}+(1-\delta\tau_k)y_{k}$.
    \STATE \;   \quad \textbf{Step 3:} Compute
    $y_{k+1}=\arg\min_{y\in \hat{K}}\{ \Phi(x_{k+1}) +\langle\nabla\Phi(x_{k+1}),y-x_{k+1}\rangle+\frac{1}{2\eta}\|y-x_{k+1}\|^2_2\}$.
    \STATE \;   \quad \textbf{Step 4:} Compute
     $z_{k+1}=\arg\min_{z\in \hat{K}}\{\alpha_k\langle \nabla\Phi(x_{k+1}),z\rangle+ \frac{1}{\delta}V_h(z,z_k)\}$.
     \STATE \; \quad  \textbf{If} A stopping condition is satisfied  \textbf{then}
      \STATE \; \quad \quad \quad\textbf{Return} $(x_k,y_k,z_k) $
       \STATE \;\quad   \textbf{end}
    \STATE  \textbf{end}
  \end{algorithmic}
\end{algorithm}
Note that the additionally updating sequence ${y_k}$ can be seen as a step of mirror descent.
The optimality condition of Algorithm 2 is given by
\begin{equation} \label{discretized.2}
\left\{
\begin{aligned}
&x_{k+1}=\delta\tau_k z_{k}+(1-\delta\tau_k)y_{k}, \\
&\|y_{k+1}-x_{k+1}\|_2=-\eta\nabla\Phi (x_{k+1}),\\
&\nabla h(z_{k+1})-\nabla h(z_{k})=-(A_{k+1}-A_k)\nabla\Phi (x_{k+1}).
\end{aligned}
\right.
\end{equation}

The following two Lemmas will be used in the sequel which plays an important role in our main results.

\begin{lemma}\label{lem:4.1}
\normalfont \cite[Sec.9]{BV}
Given an initial point $x_0\in \hat{K}$
and a sublevel set $S:=\{x\in \hat{K}:\Phi(x)\leq\Phi(x_0)\}$.
If $\Phi$ is twice continuously differentiable,
and $m_f-$strongly convex,
i.e., $\nabla^2\Phi(x)\succeq m_fI$ for all $x \in S$.
Then,
there exists an $L_{\Phi}> 0$
such that $\nabla^2\Phi(x)\preceq L_{\Phi}I$ for all $x \in S$.
\end{lemma}

\begin{lemma}\label{lem:4.2}
\normalfont
If $h$ is $\mu$-strongly convex and
$\Phi$ is $L_{\Phi}$-smooth.
Then $\|\nabla\Phi(x)\|_2\leq \frac{L_{\Phi}\sqrt{2M_{h,\mathcal{X}}}}{\sqrt{\mu}},$
where
$M_{h,\mathcal{X}} =\sup_ {x,y\in \mathcal{X}} V_h(x, y)$.
\end{lemma}
Proof. Since the optimal solution $\hat{x}^*$ of (\ref{7}) satisfies the optimality condition $\nabla\Phi(\hat{x}^* )=0$, by the smoothness of $\Phi$,
$$\|\nabla\Phi(x)\|_2\leq
\|\nabla\Phi(x)-\nabla\Phi(\hat{x}^*\|_2 +\|\nabla\Phi(\hat{x}^*)\|_2\leq
L_{\Phi}\|x-\hat{x}^*\|_2.$$
By the strong convexity of $h$, we have
$$
V_h(x,\hat{x}^*)=h(x)-h(\hat{x}^*)-\langle\triangledown h(y),x-\hat{x}^*\rangle
\geq \frac{\mu}{2}\|x-\hat{x}^*\|_2^2.$$
Combining the two inequalities above yields $\|\nabla\Phi(x)\|_2\leq
\frac{\sqrt{2}L_{\Phi}}{\sqrt{\mu}}\sqrt{V_h(x,\hat{x}^*)}.$
By Remark \ref{rem:2.1},
$M_{h,\mathcal{X}} =\sup_ {x,y\in \mathcal{X}} V_h(x, y)$.
It follows that
$\|\nabla\Phi(x)\|_2\leq \frac{L_{\Phi}\sqrt{2M_{h,\mathcal{X}}}}{\sqrt{\mu}}$.
\hfill
$\Box$

By Lemma \ref{lem:4.1},
we can assume without loss of generality that the smoothness coefficient of $\Phi$ is $L_{\Phi}$.
The following theorem shows that the convergence behavior of Algorithm 2 using a discrete-time Lyapunov function.

\begin{thm}\label{thm:4.2}
\normalfont
Suppose $h$ is $\mu$-strongly convex and
$f$ is $L_f$-smooth.
Define the discrete-time Lyapunov function
\begin{equation}
\mathcal{E}_k=A_k(\Phi(y_k)-\Phi(\hat{x}^*))+V_{h}(\hat{x}^*,z_k).
\label{LF3}
\end{equation}
Then, under Assumption \ref{a:3.1}, one has the error bound
$\frac{\mathcal{E}_{k+1}-\mathcal{E}_{k}}{\delta}\leq \varepsilon_{k+1}$,
where
$$\varepsilon_{k+1}=\frac{\delta\alpha_{k}^2}{2\mu}\|\nabla\Phi(x_{k+1})\|_2^2+\frac{A_{k+1}}{\delta}(\Phi(y_{k+1})-\Phi(x_{k+1})).$$
In addition,
if $A_k = Ck^2$, $ 0<C\leq\frac{\mu}{4L_{\Phi}}$ and $0<\eta\leq\frac{1}{kL_{\Phi}}$,
then $\varepsilon_{k+1}\leq0$
and $\|x_k-x^*\|_2=O(\frac{1}{ k}).$
\end{thm}
Proof. Using the Lyapunov function (\ref{LF3}),
we have
\begin{align*}
 \frac{\mathcal{E}_{k+1}-\mathcal{E}_{k}}{\delta}
& =-\left\langle\frac{\nabla h(z_{k+1})-\nabla h(z_{k})}{\delta},\hat{x}^*-z_{k+1}\right\rangle- \frac{1}{\delta}V_{h}(z_{k+1},z_k)+ \frac{A_{k+1}}{\delta}(\Phi(y_{k+1})-\Phi(\hat{x}^*))\\
 & \;\;\;\;-\frac{A_{k}}{\delta}(\Phi(y_k)-\Phi(\hat{x}^*))\\
& = \alpha_k\langle\nabla\Phi(x_{k+1}),\hat{x}^*-z_{k+1}\rangle-\frac{1}{\delta}V_{h}(z_{k+1},z_k) +\frac{A_{k+1}}{\delta}(\Phi(y_{k+1})-\Phi(\hat{x}^*))\\
  & \;\;\;\;
  -\frac{A_{k}}{\delta}(\Phi(y_k)-\Phi(\hat{x}^*))\\
& \leq  \alpha_k\langle\nabla\Phi(x_{k+1}),\hat{x}^*-z_{k}\rangle+ \alpha_k\langle\nabla\Phi(x_{k+1}),z_{k}-z_{k+1}\rangle -\frac{\mu}{2\delta}\|z_{k+1}-z_k\|^2 \\
  & \;\;\;\;
  +\frac{A_{k+1}}{\delta}(\Phi(y_{k+1})-\Phi(\hat{x}^*))-\frac{A_{k}}{\delta}(\Phi(y_k)-\Phi(\hat{x}^*)).
\end{align*}
The second equality follows from the third equation of (\ref{discretized.2})
and the first inequality follows from the properties of Bregman divergence,
i.e.,
$V_h(x,y)\geq \frac{\mu}{2}\|x-y\|^2_2$,
$\forall x, y\in \mathbb{R}^n$.
From Fenchel Young inequality,
we have
$\alpha_k\langle\nabla\Phi(x_{k+1}),z_{k}-z_{k+1}\rangle \leq \frac{\mu}{2\delta}\|z_{k+1}-z_k\|_2^2+\frac{\delta\alpha_k^2}{2\mu}\|\nabla\Phi(x_{k+1})\|^2_2$.
Notice that
\begin{align*}
&\frac{A_{k+1}}{\delta}(\Phi(y_{k+1})-\Phi(\hat{x}^*))-\frac{A_{k}}{\delta}(\Phi(y_k)-\Phi(\hat{x}^*))\\
=&
 \alpha_{k}(\Phi(x_{k+1})-\Phi(\hat{x}^*))+\frac{A_{k+1}}{\delta}(\Phi(y_{k+1})-\Phi(x_{k+1}))+\frac{A_{k}}{\delta}(\Phi(x_{k+1})-\Phi(y_{k})).
\end{align*}
Then it follows that
\begin{align*}
 \frac{\mathcal{E}_{k+1}-\mathcal{E}_{k}}{\delta}
& \leq  \alpha_k\langle\nabla\Phi(x_{k+1}),\hat{x}^*-z_{k}\rangle+\frac{A_{k}}{\delta}(\Phi(x_{k+1})-\Phi(y_{k}))+\alpha_k(\Phi(x_{k+1})-\Phi(\hat{x}^*))+ \varepsilon_{k+1}.
\end{align*}
From the first equation of (\ref{discretized.2}) and $A_{k+1}\leq A_{k}$, we have
\begin{align*}
 \frac{\mathcal{E}_{k+1}-\mathcal{E}_{k}}{\delta}
& \leq  \alpha_k\langle\nabla\Phi(x_{k+1}),\hat{x}^*-y_{k}\rangle+ \frac{A_{k+1}}{\delta}\langle\nabla\Phi(x_{k+1}),y_{k}-x_{k+1}\rangle+\frac{A_{k}}{\delta}(\Phi(x_{k+1})-\Phi(y_k))\\
  & \;\;\;\;
  +\alpha_k(\Phi(x_{k+1})-\Phi(\hat{x}^*))+ \varepsilon_{k+1}\\
& =  \alpha_k(\langle\nabla\Phi(x_{k+1}),\hat{x}^*-x_{k+1}\rangle+\Phi(x_{k+1})-\Phi(\hat{x}^*))+ \frac{A_{k}}{\delta}(\langle\nabla\Phi(x_{k+1}),y_{k}-x_{k+1}\rangle\\
  & \;\;\;\;
  + \Phi(x_{k+1})-\Phi(y_k))+ \varepsilon_{k+1}.
\end{align*}
This fact together with the convexity of $\Phi$ yields the conclusion $\frac{\mathcal{E}_{k+1}-\mathcal{E}_{k}}{\delta}\leq \varepsilon_{k+1}$.

We now show that $\|x_k-x^*\|_2=O(\frac{1}{k})$.
By Lemma 4.1,
there exists $L_{\Phi}> 0$
such that $\nabla^2\Phi(x)\preceq L_{\Phi}I$, $\forall x \in S$.
This upper bound on the Hessian implies for any $x_{k+1},y_{k+1}\in S$,
\begin{align}
\Phi(y_{k+1})\leq\Phi(x_{k+1})+\langle\nabla \Phi(x_{k+1}),y_{k+1}-x_{k+1}\rangle+\frac{L_{\Phi}}{2}\|y_{k+1}-x_{k+1}\|^2_2.
\label{smooth}
\end{align}
Let $0<\eta\leq\frac{1}{kL_{\Phi}}$,
plugging the second equation of (\ref{discretized.2}) into the above inequality,
we get
$\Phi(y_{k+1})-\Phi(x_{k+1})\leq -\frac{1}{2L_{\Phi}}\|\nabla \Phi(x_{k+1})\|^2_2$.
With the choices $A_k = Ck^2$ and $0 < C \leq\frac{\mu}{4L_{\Phi}}$,
the inequality
$\frac{\delta\alpha_{k}^2}{2\mu}- \frac{A_{k+1}}{2\delta L_{\Phi}}\leq 0$
holds,
this implies that $\varepsilon_{k+1}\leq0$.
Since $\Phi$ is strongly convex,
and $f$ is strongly convex and smooth,
using the
same arguments as in the second part of the proof of Theorem \ref{thm:4.1}, we can
obtain that
\begin{align}
\|y_k-x^*\|_2\leq \sqrt{\frac{2\mathcal{E}_0}{m_fA_k}}+ \frac{2}{m_f}\sqrt{\frac{2L_f(m+ \|\lambda^*\|_1)}{A_k}}.
\label{4.9}
\end{align}
Since $h$ is $\mu$-strongly convex,
by (\ref{smooth}) and Lemma \ref{lem:4.2},
$\|\nabla\Phi(x_k)\|_2\leq \frac{L_{\Phi}\sqrt{2M_{h,\mathcal{X}}}}{\sqrt{\mu}}$.
This together with the secod equation of (\ref{discretized.2}) and (\ref{4.9}) yield
\begin{align*}
\|x_k-x^*\|_2\leq &
\|x_k-y_k\|_2+\|y_k-x^*\|_2  \leq
\eta\|\Phi(x_k)\|_2+\|y_k-x^*\|_2\\ \leq &
\frac{\sqrt{2M_{h,\mathcal{X}}}}{k\sqrt{\mu}}+
\sqrt{\frac{2\mathcal{E}_0}{m_fA_k}}+ \frac{2}{m_f}\sqrt{\frac{2L_f(m+ \|\lambda^*\|_1)}{ A_k}}.
\end{align*}
Therefore, $\|x_k-x^*\|_2=O(\frac{1}{ k})$.
\hfill $\Box$

When $y_{k+1}=x_{k+1}$,
the acceleration gradient method (Algorithm 2) reduce to a gradient method,
which explicit dircretization applied to the first equation of (\ref{first order equtions}) and implicit discretization applied to the second equation of (\ref{first order equtions}).
\begin{algorithm}
  \caption{ Gradient Method (GM)}
  \label{algo:threeline}
  \begin{algorithmic}
  \STATE \;\textbf{Initialization:} Choose $x_0=z_0\in \mathbb{R}^n $, $\delta>0$.
  \STATE \;\textbf{For} $k=0,1,\cdots$ \textbf{do}
  \STATE \; \quad  \textbf{Step 1:} Set
          $A_k= \sum_{i=0}^k\frac{1}{(i+1)^2}$, $c_k=A_k$, $s_k=1/A_k$, $\alpha_k=\frac{A_{k+1}-A_{k}}{\delta}$,
           $\tau_k=\frac{A_{k+1}-A_{k}}{\delta A_{k}}$.
    \STATE \; \quad \textbf{Step 2:} Compute
    $x_{k+1}=\delta\tau_k z_{k}+(1-\delta\tau_k)x_{k}$.
    \STATE \;   \quad \textbf{Step 3:} Compute
     $z_{k+1}=\arg\min_{z\in \hat{K}}\{\alpha_k\langle \nabla\Phi(x_{k+1}),z\rangle+ \frac{1}{\delta}V_h(z,z_k)\}$.
     \STATE \; \quad  \textbf{If} A stopping condition is satisfied  \textbf{then}
      \STATE \; \quad \quad \quad\textbf{Return} $(x_k,z_k) $
       \STATE \;\quad   \textbf{end}
    \STATE  \textbf{end}
  \end{algorithmic}
\end{algorithm}
The discretization process is displayed in Algorithm 3.
The optimality condition of  Algorithm 3 is given by
\begin{equation} \label{discretized.3}
\left\{
\begin{aligned}
&x_{k+1}=\delta\tau_k z_{k}+(1-\delta\tau_k)x_{k}, \\
&\nabla h(z_{k+1})-\nabla h(z_{k})=-(A_{k+1}-A_k)\nabla\Phi (x_{k+1}).
\end{aligned}
\right.
\end{equation}

Next, we can directly give the following corollary with the proof omitted.

\begin{cor}\label{cor:4.1}
\normalfont
Suppose $h$ is $\mu$-strongly convex and $f$ is $L_f$-smooth. Then,
under Assumption \ref{a:3.1} and Lyapunov function (\ref{LF2}), one has the error bound
$\frac{\mathcal{E}_{k+1}-\mathcal{E}_{k}}{\delta}\leq \varepsilon_{k+1}$,
where
$\varepsilon_{k+1}=\frac{\delta\alpha_{k}^2}{2\mu}\|\nabla\Phi(x_{k+1})\|_2^2$.
In addition,
if $A_k = \sum_{i=0}^k\frac{1}{(i+1)^2}$,
then $\Phi(x_k)-\Phi(\hat{x}^*)\leq\frac{C_0+\mathcal{E}_0}{k}$
and the convergence rate satisfies $\|x_k-x^*\|_2=O(\frac{1}{ \sqrt{k}})$,
where $C_0=\frac{L_\Phi^2M_{h,\mathcal{X}}}{\mu^2}$,
$\mathcal{E}_0=A_0(\Phi(x_0)-\Phi(\hat{x}^*))+V_{h}(\hat{x}^*,z_0)$.
\end{cor}

\begin{remark}\label{cor:4.3}
\normalfont
If $g_i(x)=0$ for all $i\in I$,
then problem (\ref{1}) can
reduce to the problem considered in \cite{NS}.
Further,
when the full gradient is replaced by the subgradient,
Algorithm 3 is the quasi-monotone subgradient method introduced in \cite{NS}.
\end{remark}

\section{Numerical experiments}
In this section,
we solve an optimization problem with inequality constraints to illustrate the effectiveness of the accelerated gradient method  in solving such problems.
Specifically,
we conduct numerical experiments to verify the acceleration performance and effectiveness of the algorithms derived from continuous-time  dynamical system.

We consider the following quadratic optimization problem:
\begin{align}
x^*=\arg\min_{x\in\mathcal{X}} \frac{1}{2}x_1^2+\frac{3}{2}(x_2-1)^2, \;\; \mbox{s.t.}\;x_2-x_1-1\leq 0,
\label{5.1}
\end{align}
Next,
we examine two distance-generating functions $h(\cdot)$
and their corresponding constraint sets $\mathcal{X}$:
\begin{itemize}
\item Squared Euclidean norm: $h(x)=\frac{1}{2}\|x\|_2^2$ and $\mathcal{X}=\{x:\|x\|_2\leq R\}$ for some constant $R>0$.
\item Negative entropy: $h(x)=\sum_{i=1}^n x_i\log x_i$ and $\mathcal{X}=\{x:\sum_{i=1}^n x_i=1,x_i\geq 0\; \mbox{for}\; i=1,\cdots,n\}$.
\end{itemize}
In the following simulation,
we show how to track the optimal solution $x^*$ using Algorithm 2.
The augmented objective function (\ref{8}) takes the form
$$\Phi(x):= \frac{1}{2}x_1^2+\frac{3}{2}(x_2-1)^2- \frac{1}{c }\log(s +x_1-x_2+1).$$
For the squared Euclidean norm distance,
the radius of the Euclidean ball $\mathcal{X}$ is set to be $R=3$,
the initial points $z$ and $y$ are randomly generated in $\mathcal{X}$,
$\eta=1/4$
and
$C=1/10$.
For the negative entropy distance,
the initial points $z$ and $y$ are randomly generated in a $2$-dimensional positive simplex,
$\eta=1/2$
and
$C=2/3$.
Moreover,
using particular values,
all conditions of Theorem \ref{thm:4.2} are satisfied.

\begin{figure}[H]
    \centering
    \subfigure{
    \includegraphics[width=0.5\textwidth]{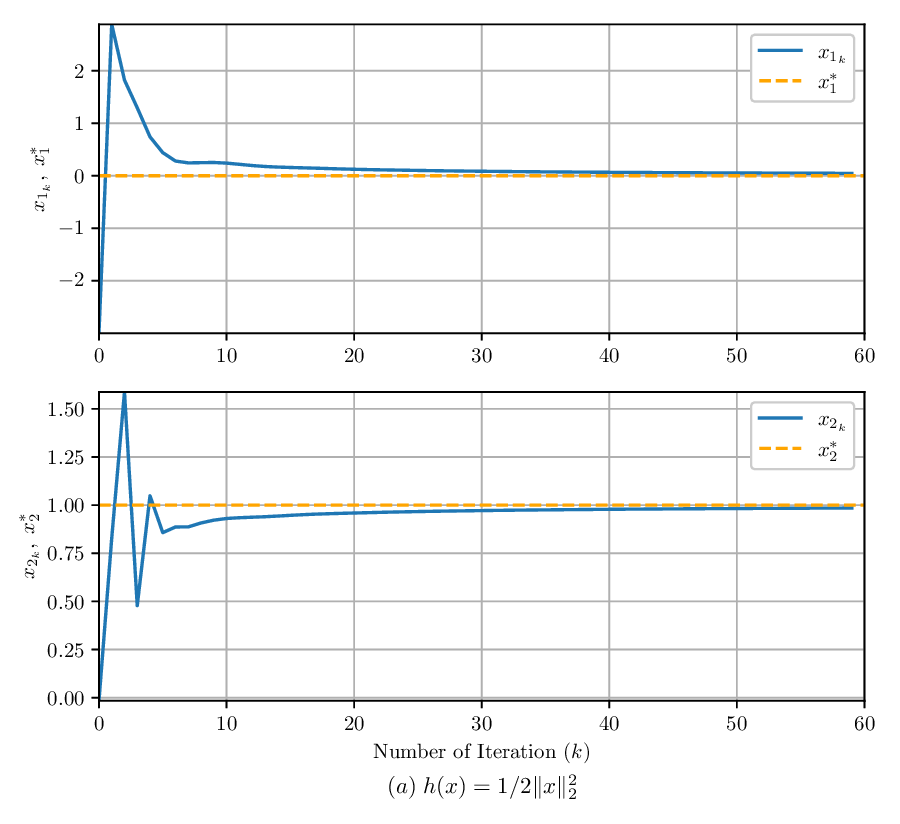}
    \includegraphics[width=0.5\textwidth]{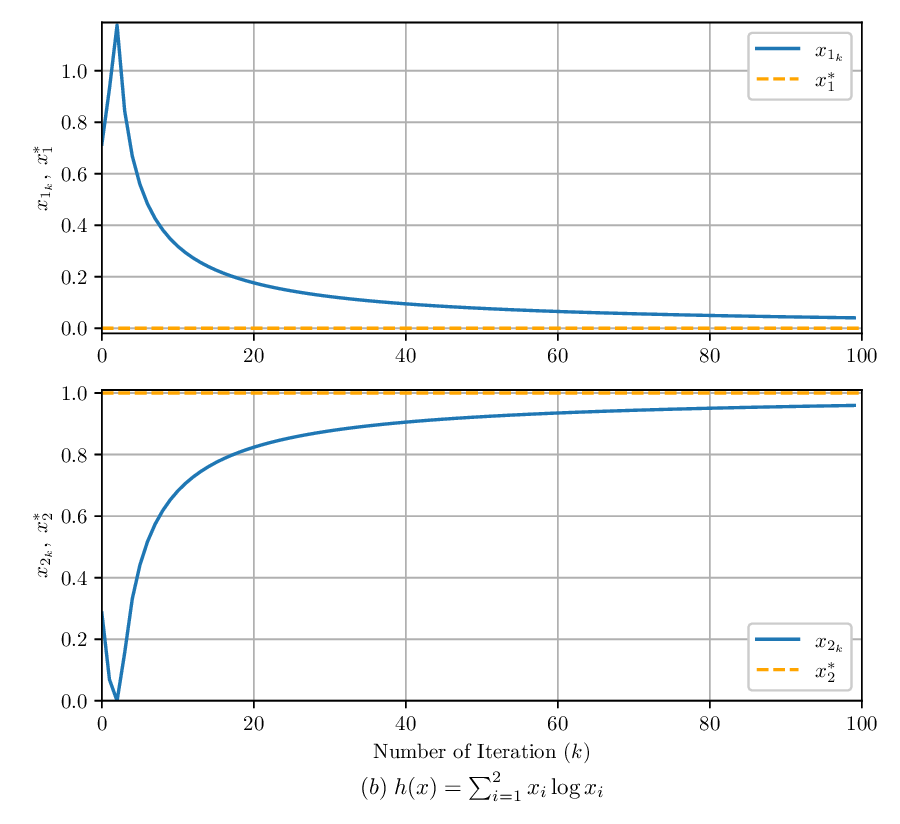}}
    \caption{the coordinates of the optimal trajectory $x^*$ defined in (5.1) and the tracking trajectory $x_k$.}
    \label{fig:Figure 5.1}
\end{figure}

\begin{figure}[ht]
   \centering
    \subfigure{
    \includegraphics[width=0.5\textwidth]{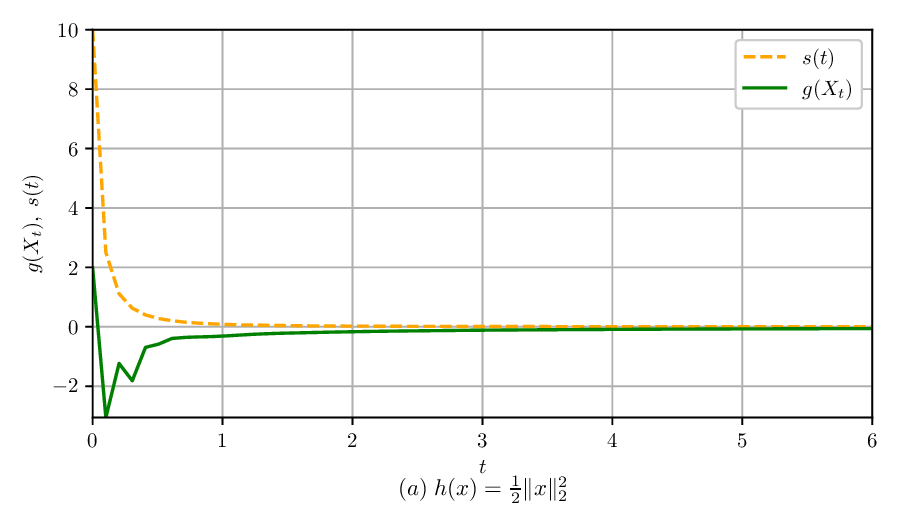}
    \includegraphics[width=0.5\textwidth]{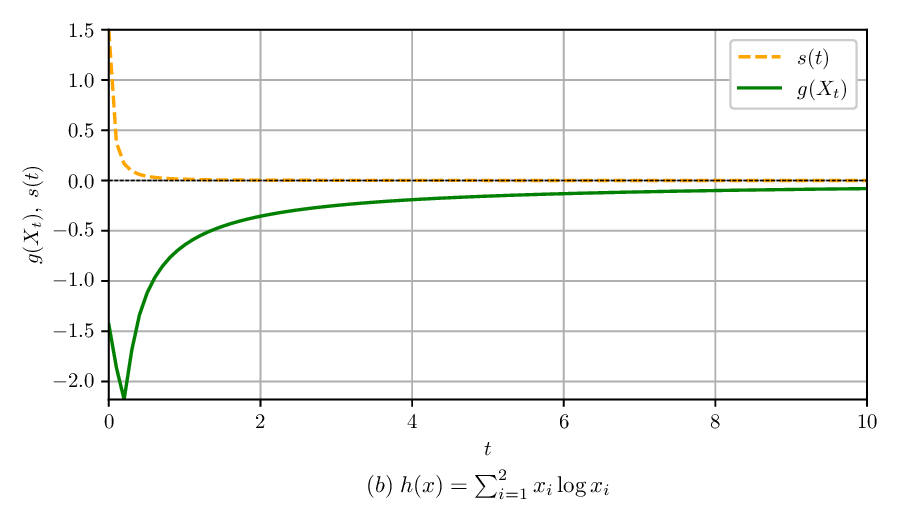}}
    \caption{the coordinates of the constrain function $g(X_t)$ and the slack variable $s(t)$ against t, for the problem (5.1).}
    \label{fig:Figure 5.2}
\end{figure}

\begin{figure}[ht]
\centering
    \subfigure{
    \includegraphics[width=0.5\textwidth]{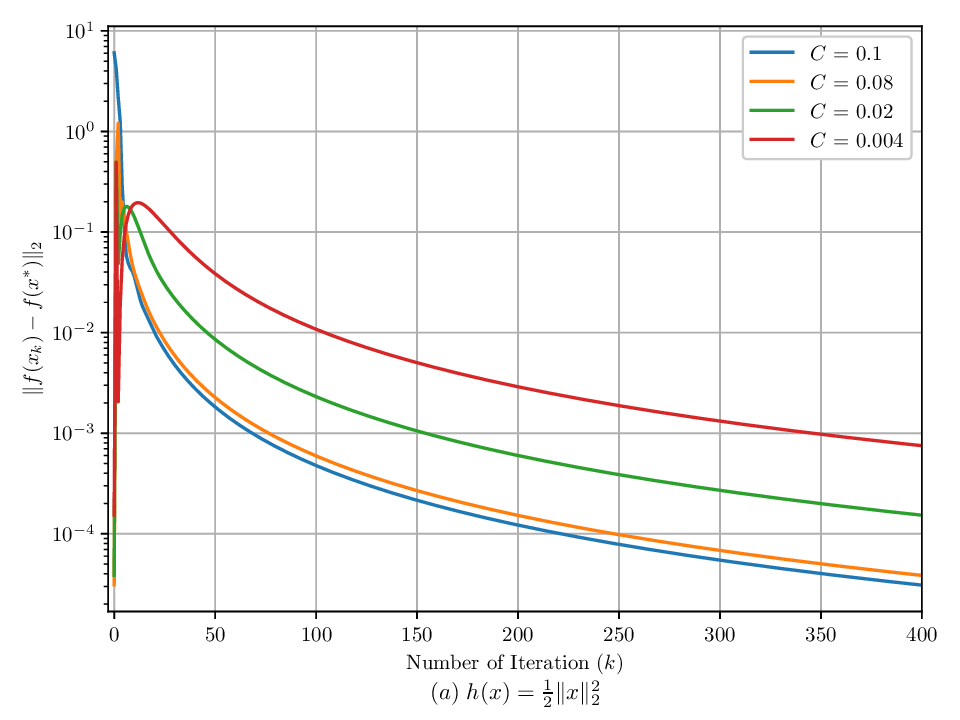}
    \includegraphics[width=0.5\textwidth]{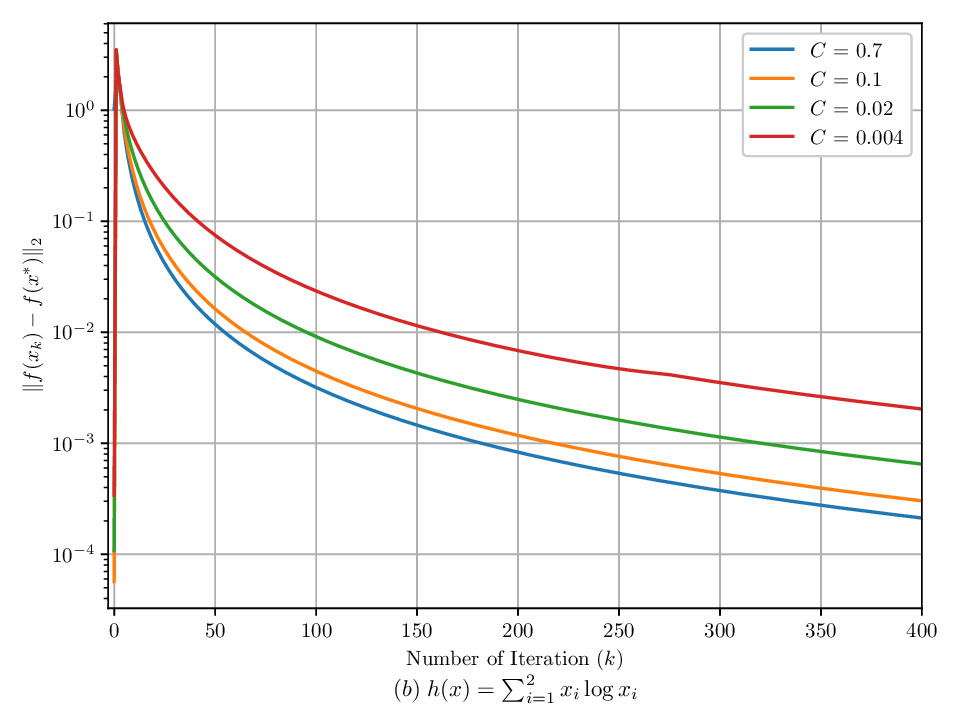}}
    \caption{the performance of Algorithm 2 applied to (\ref{5.1}) under different distance generating functions.}
    \label{fig:Figure 5.3}
\end{figure}

The trajectory of the solution $x_k=(x_{1_{k}},x_{2_{k}})^T$, along with the optimal solution $x^*=(x_1^*,x_2^*)$ defined in (\ref{5.1}), are shown in Figure \ref{fig:Figure 5.1}.
In detail,
Figure \ref{fig:Figure 5.1}(a) and Figure \ref{fig:Figure 5.1}(b) show the results of the squared Euclidean norm distance and the negative entropy distance for the quadratic optimization problem,
respectively.
Figure \ref{fig:Figure 5.2}(a) and Figure \ref{fig:Figure 5.2}(b) display the time evolution of the constraint function $g(x):=x_2-x_1-1$ and the nonnegative time-dependent slack function $s(t)$ about the squared Euclidean norm distance and the negative entropy distance,
respectively.
We can see the state $X_t$ violates the constraint $g(x)\leq 0$ at $t =0$.
However,
$X_t$ converges to the feasible set exponentially fast as the slack function $s(t)$ vanishes exponentially.
In Figure \ref{fig:Figure 5.3},
we illustrate the performance of Algorithm 2 applied to (\ref{5.1}), across various parameter settings for $C$.
We observe that the large $C$ is,
the better  Algorithm 2 performs.

\begin{figure}[ht]
\centering
\includegraphics[width=.6\textwidth]{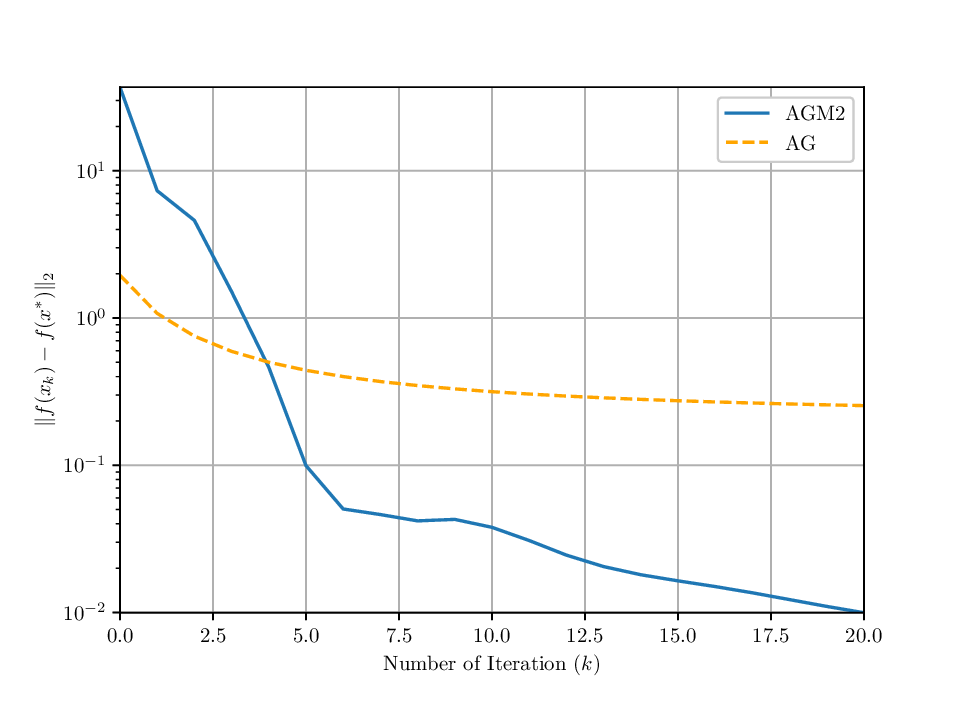}
    \caption{The objective function value gap obtained by Algorithm 2 (AGM2) and Algorithm 3 (AG).}
    \label{fig:Figure 5.4}
\end{figure}

Finally,
we consider
$\eta=1/4$ and
$C=1/10$
with the squared Euclidean norm distance.
The radius of the Euclidean ball $\mathcal{X}$ is set to be $R=1$.
For the stopping criterion,
we choose $\|f(x_k)-f(x^*)\|_2 \leq 10^{-2}$.
In Figure \ref{fig:Figure 5.4},
we compare the performance of  Algorithm 2 (AGM2) against  Algorithm 3 (AG).
We can clearly observe the effect that acceleration has on reducing the number of required iterations to obtain a certain accuracy.
We noticed that,
during the application of the Quasi-monotone method,
the selection of $A_k$ with a relatively small magnitude can lead to a series of issues.
Specifically,
this small magnitude may lead to significant increases in rounding errors and cumulative errors,
thereby weakening the stability of numerical calculations and resulting in oscillations characterized by repeated increases and decreases in error values.

\section{Conclusions}
In this paper,
we investigated the accelerated gradient method for convex optimization problems with inequality constraints.
Specifically,
we approximate the inequality constrained convex optimization problem as an unconstrained optimization problem by employing the logarithmic barrier function.
Under ideal scalar assumptions,
we proposed a continuous-time dynamical system associated with the Bregman Lagrangian
and derived its polynomial and exponential convergence rates based on the relevant parameter settings.
Additionally,
we developed several discrete-time algorithms by discretizing the continuous-time dynamical system,
and obtained the fast convergence results matching that the underlying dynamical system.
The numerical experiments demonstrated the acceleration performance and effectiveness of the proposed algorithms.

\end{document}